\documentclass[12pt]{article}

\usepackage{amsmath}
\usepackage{amssymb}
\usepackage{amsthm}

\usepackage[dvips]{graphicx}
\usepackage{color}
\usepackage{epsfig}
\usepackage[mathscr]{eucal}
\newtheorem{theorem}{Theorem}

\newtheorem{prop}[theorem]{Proposition}

\newtheorem{lemma}[theorem]{Lemma}
\newtheorem{cor}[theorem]{Corollary}

\newcommand{\wt}{\ensuremath{\textrm{wt}}}
\newcommand{\ft}{\ensuremath{\mathbb{F}_q}}
\newcommand{\ftn}{\ensuremath{\mathbb{F}_{q^n}}}
\newcommand{\pr}{\ensuremath{\textrm{Pr}}}

\allowdisplaybreaks

\title{De Bruijn Cycles for Covering Codes}
\author{Fan Chung and Joshua N. Cooper \\ \small Department of Mathematics \\ \small University of California, San Diego, La Jolla, CA}
\date{\today}

\begin{document}

\maketitle

\begin{abstract}
A de Bruijn covering code is a $q$-ary string $S$ so that every $q$-ary string is at most $R$ symbol changes from some $n$-word appearing consecutively in $S$.  We introduce these codes and prove that they can have length close to the smallest possible covering code.  The proof employs tools from field theory, probability, and linear algebra.  We also prove a number of ``spectral'' results on de Bruijn covering codes.  Included is a table of the best known bounds on the lengths of small binary de Bruijn covering codes, up to $R=11$ and $n=13$, followed by several open questions in this area.
\end{abstract}

\section{Introduction} \label{section:introduction}

A covering code $\mathcal{C}$ of radius $R$ and dimension $n$ on $q$ symbols is a subset of the space $[q]^n$ such that every string in $[q]^n$ differs from some element of $\mathcal{C}$ in at most $R$ coordinates.  It is common to require that $R$ be as small as possible in the definition of a covering code, but, for the sake of notational convenience, we do not require this here.\\

\noindent Question: Given $n$, $R$, and $q$, what is the smallest $M=M(n,R,q)$ so that there exists an $q$-ary string $S=(s_0,\ldots,s_{M-1})$ with the property that the set of $n$-strings appearing as $(s_i,\ldots,s_{i+n-1})$, with indices taken modulo $M$, form a covering code of radius $R$?  Call such a string a $(n,R,q)$-{\it de Bruijn covering code}.

For example, $111000$ is a $(4,1,2)$-de Bruijn covering code, because every binary $4$-string is at most one bit change from an element of
$$\{1110,1100,1000,0001,0011,0111\}.$$
On the alphabet $\{\textrm{A,G,T,C}\}$, the string
$$
\textrm{AGATCGCAGATATGGTCTATG}
$$
is a $(4,2,4)$-de Bruijn covering code, by Proposition \ref{enlarge} below.

Clearly, $M(n,0,q) = q^n$, since any de Bruijn covering code of radius $0$ is actually a de Bruijn cycle, and de Bruijn cycles of all orders over an arbitrary alphabet exist.  (See, for example, \cite{Land00}.)  If we fix $R>0$ and $q \geq 2$, how does $M(n,R,q)$ grow as $n \rightarrow \infty$?

It is easy to see that the growth is at least $\Omega(q^n/n^R)$, by the so-called ``sphere-covering'' bound.  The set of strings which differ from any given $S$ in at most $R$ places has the same cardinality, $\sum_{k=0}^R \binom{n}{k} (q-1)^k$.  Therefore, if we are to cover all $q^n$ strings, we need at least
$$
\frac{q^n}{\sum_{k=0}^R \binom{n}{k} (q-1)^k }
$$
codewords.  On the other hand, it is well known that the size of the smallest $q$-ary covering code of radius $R$ actually achieves this bound, up to a multiplicative constant which depends on $R$ and $q$.  (See \cite{KSV03} for the latest results on the size of this constant.)  We may concatenate all the codewords of such a minimal code to yield a $(n,R,q)$-de Bruijn covering code of length $O(q^n/n^{R-1})$.  This construction is clearly very wasteful, however.  Can we do better, i.e., is the true order of magnitude of $M(n,R,q)$ closer to the sphere-covering bound?  In particular, can we say something nontrivial in the case of $R=1$?  In fact, in Section \ref{section:mainresult} we prove the following.

\begin{theorem} \label{mainthm} For each $n$ and $q$ a prime power, there exists a $(n,R,q)$-de Bruijn covering code of length $\leq (R+1 + o(1)) q^n \log n  / (\binom{n}{R} (q-1)^{R})$.
\end{theorem}

Section \ref{section:preliminaries} states several definitions and preliminary results we will need to prove this.  The next section contains the proof itself, and Section \ref{section:spectral} introduces a ``spectral'' perspective on de Bruijn covering codes that holds some independent interest.  In Section \ref{section:bounds} we present bounds for special values of $n$, $R$, and $q$, and include a table of bounds on $M(n,R,2)$ for $2 \leq n \leq 13$ and $1 \leq R \leq 11$.  We end with several remarks and questions for further work in Section \ref{section:remarksandfurtherquestions}.

\section{Preliminaries} \label{section:preliminaries}

We fix a prime power $q \geq 2$ throughout this section and the next, and take our alphabet to be $\ft$.  (If $q$ is not a prime power, we take the alphabet to be $\mathbb{Z}/q\mathbb{Z}$.) Write $b_R(v)$, for $v$ an $n$-string drawn from $\ft$, to denote the set of those strings differing from $v$ in at most $R$ coordinates.  That is, $b_R(v)$ is $v$'s radius $R$ neighborhood in the Hamming metric.  Also, write $\wt(v)$ for the Hamming weight of the vector $v$, the number of nonzero symbols it contains.

Let $\alpha$ be a generator of the multiplicative group of the finite field $\ftn$.  Denote by $\mathcal{E}$ the elementary basis for $\ft^n$ over $\ft$.  Given a basis $\mathcal{B} = \{b_1,\ldots,b_n\}$ of $\ftn$ over $\ft$ and an element $\gamma \in \ftn$, write $f_{\mathcal{B}}(\gamma)$ for the element of $\ft^n$ whose $j^\textrm{th}$ coordinate is the coefficient of $b_j$ in the $\mathcal{B}$-representation of $\gamma$.  Then, given a nonzero vector $\mathbf{x} \in \ft^n$, define $\Lambda(\alpha,\mathcal{B},\mathbf{x})$ to be the string whose $j^\textrm{th}$ coordinate (i.e., $\Lambda_j(\alpha,\mathcal{B},\mathbf{x})$, $1 \leq j \leq q^n-1$) is $\mathbf{x}^\intercal f_{\mathcal{B}}(\alpha^j)$.  
It is well known that, when $\mathcal{B}=\{\alpha^j : 0 \leq j \leq n-1\}$ and $\wt(\mathbf{x})=1$, $\Lambda(\alpha,\mathcal{B},\mathbf{x})$ is a de Bruijn cycle of order $n$ if we insert a $0$ at the beginning.  (See, for example, \cite{Fred82}.)  We generalize this result as follows.  Define $\Lambda^*(\alpha,\mathcal{B},\mathbf{x})$ to be the sequence $\Lambda(\alpha,\mathcal{B},\mathbf{x})$ with a zero inserted at the beginning of each occurrence of the string $0\ldots01$.  Then we have the following.

\begin{prop} \label{hom} Fix a basis $\mathcal{B}$ of $\ft^n$ over $\ft$, a generator $\alpha \in \ftn^\times$, and a vector $\mathbf{x} \in \ft^n$, and write $\Phi(j)$ for the vector
$$
(\Lambda_j(\alpha,\mathcal{B},\mathbf{x}),\ldots,\Lambda_{j+n-1}(\alpha,\mathcal{B},\mathbf{x}))^\intercal \in \ft^n
$$
The map $\Psi$ which sends $0$ to $0$ and $\alpha^j$ to $\Phi(j)$ is an isomorphism from the additive group of $\ftn$ to $\ft^n$.
\end{prop}
\begin{proof} First, we show that $\Psi$ is linear.  Write $e_j$ for the elementary $n$-vector whose coordinates are all zero except for a $1$ in the $j^\textrm{th}$ coordinate.  We denote by $M_{\gamma,\mathcal{B}}$ the matrix representing multiplication by $\gamma \in \ftn$ in the $\mathcal{B}$ basis.  It is easy to see that
$$
\Lambda_j(\alpha,\mathcal{B},\mathbf{x}) = \mathbf{x}^\intercal f_{\mathcal{B}}(\alpha^j)
$$
and therefore that
\begin{equation} \label{matrixformula}
\Psi(\gamma) = \sum_{j=0}^{n-1} e_{j+1} \mathbf{x}^\intercal f_{\mathcal{B}}(\alpha^j \gamma) = \sum_{j=0}^{n-1} e_{j+1} \mathbf{x}^\intercal M_{\alpha,\mathcal{B}}^j f_{\mathcal{B}}(\gamma),
\end{equation}
which is obviously linear.

Now, suppose that $\Psi(\gamma)=0$.  We show that $\gamma=0$.  Indeed, suppose that $\{j_1,\ldots,j_n\}$ are $n$ distinct integers so that $\Lambda_{j_i}(\alpha,\mathcal{B},\mathbf{x}) = 0$ for each $i$.  If we denote by $S$ the subspace of $\ft^n$ orthogonal to $\mathbf{x}$, then we have $\alpha^{j_i} \in f_{\mathcal{B}}^{-1}(S)$ for each $i$.  However, $f_{\mathcal{B}}$ is linear and has a trivial kernel, so all the $\alpha^{j_i}$ lie in a subspace of $\ft^n$ of dimension $n-1$ and are therefore linearly dependent.  If we take $j_i=j+i$ for some $j$ (i.e., $\Psi(\gamma)=0$ with $\gamma = \alpha^{j}$), then we have that $\{\alpha^i\}_{i=j+1}^{j+n}$ is a dependent set.  Since $M_{\alpha,\mathcal{B}}$ is nonsingular, this implies that $\{\alpha^i\}_{i=0}^{n-1}$ is a dependent set.  But then we have
$$
\sum_{i=0}^{n-1} c_i \alpha^i = 0
$$
for some nonzero $(c_1,\ldots,c_n)$, so $\alpha$ satisfies a polynomial identity of degree less than $n$.  Since $\alpha$ generates $\ftn^\times$, this implies that $\{\alpha^j\}_{j=0}^{d}$ is a basis for $\ftn$ for some $d < n-1$, contradicting the fact that the dimension of $\ftn$ over $\ft$ is $n$.  We can therefore conclude that $\gamma=0$.
\end{proof}

Note that the map $\gamma \mapsto M_{\gamma,\mathcal{B}}$ is actually an isomorphism of fields.  The image is a set of matrices which form a field, i.e., a \textit{matrix field}.  These objects have been studied extensively and thoroughly characterized when the matrices take their entries from a finite field (\cite{Beard86}).

\begin{cor} $\Lambda^*(\alpha,\mathcal{B},\mathbf{x})$ is a de Bruijn cycle.
\end{cor}
\begin{proof} By the above argument, $\Lambda(\alpha,\mathcal{B},\mathbf{x})$ contains all nonzero $n$-strings.  Clearly, the insertion of a $0$ causes the occurrence of the all-zeroes string without disrupting the presence of any other string.
\end{proof}

Our approach is to find an $\alpha \in \ftn$, a basis $\mathcal{B}$, and a vector $\mathbf{x}$ so that the first $K \sim q^n \log n / (\binom{n}{R} (q-1)^R)$ length $n$ strings appearing in $\Lambda(\alpha,\mathcal{B},\mathbf{x})$ are (almost) a covering code of radius $R$.  Specifically, we wish to show that, for only a small fraction of all $v \in \ftn$,
$$
(v + B_R(0^n)) \cap \Psi(\{\alpha^j\}_{j=1}^K) = \emptyset
$$
where $\Psi$ is the function defined in Proposition \ref{hom}.  Define $\Psi^\prime = f_\mathcal{B} \circ \Psi^{-1}$.  Setting $w = \Psi^{-1}(v)$, we may bound this quantity from above by asking the number of $w$ so that
$$
\Psi^\prime \left(\binom{\mathcal{E}}{R} \right) \cap f_\mathcal{B}(w + \{\alpha^j\}_{j=1}^K) = \emptyset
$$
which, by (\ref{matrixformula}), is the same as saying that 
$$
\left \{ \left ( \sum_{j=0}^{n-1} e_{j+1} \mathbf{x}^\intercal M_{\alpha,\mathcal{B}}^j \right )^{\!-1}\!\!\!\!v : \wt(v)=R \right \} \cap f_\mathcal{B}(w + \{\alpha^j\}_{j=1}^K) = \emptyset.
$$
We must determine which matrices may appear in the form of the left-hand term.  First, a result from linear algebra is needed.  The following theorem appears in \cite{Beard86}.  A \textit{non-derogatory} matrix is one whose eigenspaces are all one-dimensional, and a matrix in \textit{rational canonical form} is comprised of blocks of the form

$$
\begin{array}{cccccc}
0 & 0 & 0 & 0 & \cdots & a_1 \\
1 & 0 & 0 & 0 & \cdots & \vdots \\
0 & 1 & 0 & 0 & \cdots & \vdots \\
0 & 0 & 1 & 0 & \cdots & \vdots \\
\vdots & \vdots & \vdots & \vdots & \ddots & \vdots \\
0 & \cdots & \cdots & \cdots & 1 & a_n \\
\end{array}
$$
along the diagonal.

\begin{theorem} \label{ratcan} If $A \in K^{n \times n}$ is non-derogatory and in rational canonical form, then the following are equivalent:
\begin{enumerate}
\item $X$ commutes with $A$.
\item The successive columns of $X$ are $v$, $Av$, \ldots, $A^{n-1}v$ for any $v \in K^n$.
\item There exists a polynomial $g \in K[x]$ so that $X = g(A)$.
\end{enumerate}
Furthermore, $g = \sum_{j=0}^{n-1} v_{j+1} x^j$.
\end{theorem}

The matrices $M_{\alpha,\mathcal{B}}$ are non-derogatory when $\alpha$ is a generator of $\ftn$, because their eigenvalues are all distinct, as the next result states.

\begin{prop} \label{eigen} A matrix $M \in \ft^{n \times n}$ is of the form $M_{\alpha,\mathcal{B}}$ for some generator $\alpha \in \ftn^\times$ and basis $\mathcal{B} \subset \ftn$ over $\ft$ if and only if its eigenvalues (over the algebraic closure of $\ft$) are $\{\alpha^{q^j}\}_{j=0}^{n-1}$.
\end{prop}
\begin{proof} For a given $\alpha$, fix the basis $\mathcal{A} = \{\alpha^j\}_{j=0}^{n-1}$.  Clearly, if we write $B$ for the matrix whose columns are $\mathcal{B}$ written in the basis $\mathcal{A}$, then $M_{\alpha,\mathcal{B}} = B^{-1} M_{\alpha,\mathcal{A}} B$.  Therefore, a matrix $M$ is one of the desired ones if and only if it has the same eigenvalues as the matrix $M_{\alpha,\mathcal{A}}$.  Let $p_\alpha(\lambda)$ denote the characteristic polynomial of this matrix.  By the Cayley-Hamilton Theorem (which applies to all commutative rings), $p_\alpha(M_{\alpha,\mathcal{A}}) = 0$.  However, the map $\alpha \mapsto M_{\alpha,\mathcal{B}}$ is an isomorphism of fields for any basis $\mathcal{B}$.  Therefore, $p_\alpha(\alpha)=0$.  Since the Galois group of $\ftn$ over $\ft$ is cyclic and generated by the Frobenius map $x \mapsto x^q$, and the rest of the roots of $p_\alpha$ are the Galois conjugates of $\alpha$, the result follows.
\end{proof}

Furthermore, if we let $\Theta_\alpha$ denote the basis $\{\alpha^j\}_{j=0}^{n-1}$, then $M_{\alpha,\Theta_\alpha}$ is in rational canonical form.  Its $j^\textrm{th}$ column is $e_{j+1}$ for $1 \leq j \leq n-1$ and its $n^\textrm{th}$ column is the vector of coefficients of the minimal polynomial of $\alpha$ (without the leading term).  Using this fact, we can prove the following from Theorem \ref{ratcan}.

\begin{lemma} \label{unif} Fix a generator $\alpha$ of $\ftn$.  Choose $\mathbf{x} \in \ft^n \setminus \{0^n\}$ randomly and uniformly, and choose a basis $\mathcal{B}$ randomly and uniformly.  Then
$$
\left ( \sum_{j=0}^{n-1} e_{j+1} \mathbf{x}^\intercal M_{\alpha,\mathcal{B}}^j \right )^{\!-1}
$$
is distributed uniformly over all invertible matrices.
\end{lemma}
\begin{proof} Evidently, it suffices to show that $D(\mathcal{B},\mathbf{x}) = \sum_{j=0}^{n-1} e_{j+1} \mathbf{x}^\intercal M_{\alpha,\mathcal{B}}^j$ is distributed uniformly.  This matrix is one whose rows are $\mathbf{x}^\intercal$, $\mathbf{x}^\intercal M_{\alpha,\mathcal{B}}$, $\ldots$, $\mathbf{x}^\intercal M_{\alpha,\mathcal{B}}^{n-1}$.  Write $A$ for the matrix $M_{\alpha,\Theta_\alpha}$ and $P$ for the matrix whose successive columns are the elements of $\Theta_\alpha$ written in the $\mathcal{B}$ basis, and write $\mathbf{y}$ for $P^\intercal \mathbf{x}$.  Then we may also say that $D(\mathcal{B},\mathbf{x})$ is the matrix whose rows are $\mathbf{x}^\intercal$, $\mathbf{x}^\intercal P A P^{-1}$, $\ldots$, $\mathbf{x}^\intercal P A^{n-1} P^{-1}$, which we may rewrite as $D(A,P^\intercal \mathbf{x}) P^{-1}$.    Therefore, by Theorem \ref{ratcan} and the fact that $A$ is non-derogatory and in rational canonical form, $D(\mathcal{B},\mathbf{x}) = g_{y}(A)^\intercal P^{-1}$ with $g_\mathbf{y}$ denoting the polynomial whose coefficients are the entries of $\mathbf{y}$.  Choosing $\mathbf{x}$ uniformly and randomly from the nonzero vectors yields the same distribution on $\mathbf{y}$, independent of the choice of $\mathcal{B}$.  Since $A$ is the image of $\alpha$ under the map $\alpha \mapsto M_{\alpha,\Theta_\alpha}$, and $g_y(\alpha)$ is uniformly distributed over $\ftn \setminus \{0\} $ as $\mathbf{y}$ varies, we have $g_y(A)$ uniformly distributed over all matrices of the form $M_{\gamma,\Theta_\alpha}$ for $\gamma \in \ftn \setminus \{0\}$.  Choosing $\mathcal{B}$ uniformly is the same as choosing $P^{-1}$ uniformly, so we may conclude that $D(\mathcal{B},\mathbf{x}) = g_{y}(A)^\intercal P^{-1}$ is uniformly distributed over all invertible matrices.

\end{proof}

\section{The Main Result} \label{section:mainresult}

It remains to show that the set of all sums of $k$ columns of a randomly, uniformly chosen invertible matrix are distributed more or less uniformly.  Before proceeding, we need to state Suen's Inequality.  We follow \cite{AS00}.  Let $\{A_i\}_{i \in I}$ be a set of events, and define a symmetric relation (i.e, a graph) $\sim$ on $I$.  We say that $\sim$ is a \textit{superdependency} graph if, whenever $J_1, J_2 \subset I$ have no edges between them, any Boolean combination of $\{A_i\}_{i \in J_1}$ is independent of any Boolean combination of $\{A_i\}_{i \in J_2}$.  Write $M = \prod_{i \in I} \pr[\overline{A_i}]$.

\begin{theorem}[Suen's Inequality] Define
$$
y(i,j) = (\pr[A_i \wedge A_i] + \pr[A_i]\pr[A_j]) \prod_{l \sim i \textit{ or } l \sim j} (1 - \pr[\overline{A_l}])^{-1}.
$$
Then
$$
\pr \left [ \bigwedge_{i \in I} \overline{A_i} \right ] \leq M e^{\sum_{i \sim j} y(i,j)}.
$$
\end{theorem}

The following is a routine application of this result.

\begin{prop} \label{invertible} For $R \in \mathbb{Z}^+$, if $M$ is chosen randomly and uniformly from $GL_n(\ft)$, then, for any set $S \subset \ft^n$ with $|S| = q^n K/(\binom{n}{R} (q-1)^R)$,
$$
\pr \left[\left \{ M v : \wt(v)=R \right \} \cap S = \emptyset \right ] \leq e^{-K} (c_q^{-1}+o(1)).
$$
where $c_q = \prod_{j=1}^\infty (1-q^{-j})$ and $K = o(\sqrt{n})$.
\end{prop}
\begin{proof} The probability that a randomly, uniformly chosen invertible matrix has all sums of $k$ columns lying outside of a set $S$ is given by
\begin{align*}
\rho &= \pr [M v \in \overline{S} \textrm{ when } \wt(v)=R | M \in GL_n(\ft)] \\
& = \frac{\pr [(M v \in \overline{S} \textrm{ when } \wt(v)=R) \wedge (M \in GL_n(\ft))]}
{\pr[M \in GL_n(\ft)]} \\
& \leq \frac{\pr [M v \in \overline{S} \textrm{ when } \wt(v)=R]}
{\pr[M \in GL_n(\ft)]}
\end{align*}
where we are choosing $M$ randomly and uniformly from \textit{all} matrices.  It is well known that $|GL_n(\ft)| = q^{n^2} (c_q+o(1))$ with $c_q = \prod_{j=1}^\infty (1-q^{-j})$.  Therefore,
$$
\rho \leq \pr [M v \in \overline{S} \textrm{ when } \wt(v)=R] (c_q^{-1} + o(1)).
$$
Now, for a vector $v$ of weight $R$, define $A_v$ to be the event that $M v \in S$, and let $I(v)$ denote the set of indices at which $v$ is nonzero.  Then $\pr[M v \in \overline{S} \textrm{ when } \wt(v)=R] = \pr[\wedge_{v} \overline{A_v}]$.  The relation $v \sim w$ iff $I(v) \cap I(w) \neq \emptyset$ clearly defines a superdependency graph on these events.  Furthermore, any pair $A_{v}$ and $A_{w}$, $v \neq w$, are independent, since, if we fix the $i^\textrm{th}$ columns of $M$ for $i \in I(v) \cap I(w)$, then $\sum_{i \in I(v) \setminus I(w)} M e_i$ and $\sum_{i \in I(v) \setminus I(w)} M e_i$ are independent and uniformly distributed over $\ft^n$.  Therefore,
\begin{align*}
y(v,w) &= 2 \,\pr[A_{v}]\pr[A_{w}] \prod_{z \sim v \textrm{ or } z \sim w} (1 - \pr[\overline{A_{z}}])^{-1} \\
&\leq 2 \left (\frac{K}{\binom{n}{R}(q-1)^R} \right )^2 \left (1 - \frac{K}{\binom{n}{R}(q-1)^R} \right)^{-2 \left (\binom{n}{R}(q-1)^R - \binom{n-R}{R}(q-1)^R \right)} \\
&= 2 \left (\frac{K}{\binom{n}{R} (q-1)^R} \right )^2 \left (1 - \frac{K}{\binom{n}{R}(q-1)^R} \right)^{\binom{n}{R-1} (-2R^2+o(1))(q-1)^R} \\
&\leq 2 \left (\frac{K}{\binom{n}{R}(q-1)^R } \right )^2 e^{ - K \binom{n}{R-1} (-2R^2+o(1))/\binom{n}{R}} \\
&= 2 \left (\frac{K}{\binom{n}{R}(q-1)^R}  \right )^2 e^{ - K (-2R^3+o(1))/n}.
\end{align*}
Since there are $\binom{n}{R} \left(\binom{n}{R} - \binom{n-R}{R}\right)(q-1)^{2R}/2 = O(n^{2R-1})$ relations $v \sim w$, the quantity $\sum_{v \sim w} y(v,w)$ tends to $0$ as $n \rightarrow \infty$ so long as $K = o(\sqrt{n})$.  Therefore, Suen's Inequality implies that
\begin{align*}
\pr \left [ \bigwedge_{\wt(v)=R} \overline{A_v} \right ] & \leq  (c_q^{-1}+o(1)) \prod_{\wt(v)=R} \pr[\overline{A_v}]\\
& = (c_q^{-1}+o(1)) \left(1 - \frac{K}{\binom{n}{R} (q-1)^R}\right)^{\binom{n}{R} (q-1)^R} \\
& \leq (c_q^{-1}+o(1)) e^{-K}.
\end{align*}
\end{proof}

Taking an initial segment of a random $\Lambda(\alpha,\mathcal{B},\mathbf{x})$ and adding in all the ``uncovered'' codewords yields an $(n,R,q)$-de Bruijn covering code.

\setcounter{theorem}{1}
\begin{theorem} For each $n$, there exists an $(n,R,q)$-de Bruijn covering code of length $\leq (R+1 + o(1)) q^n \log n / (\binom{n}{R} (q-1)^R)$.
\end{theorem}
\begin{proof} Fix any generator $\alpha \in \ftn^\times$.  Choose the basis $\mathcal{B}=\{b_i\}_{i=1}^n$ and the vector $\mathbf{x} \in \ft^n \setminus \{0^n\}$ randomly and uniformly.  Then define $\overline{\Lambda}(K)$ to be the string of the first $q^n K / (\binom{n}{R} (q-1)^R) + n$ symbols of $\Lambda(\alpha,\mathcal{B},\mathbf{x})$ (which we will call $\Lambda_1(K)$), followed by a concatenated list (which we will call $\Lambda_2(K)$) of all strings in
$$
\ft^n \setminus \bigcup_{c \in \mathcal{C}} b_R(c)
$$
where $\mathcal{C}$ is the set of codewords appearing as $n$ consecutive symbols (without wrap-around) in $\Lambda_1(K)$.  Then the resulting expected length of the string is given by
\begin{equation} \label{eq2}
\textrm{E}(|\Lambda_1(K)|+|\Lambda_2(K)|) = \frac{q^n K}{\binom{n}{R} (q-1)^R} + n + n q^n \sum_{v \in \ft^n} \pr[b_R(v) \cap \mathcal{C} = \emptyset]
\end{equation}
Furthermore, the constructed string is an $(n,R,q)$-de Bruijn covering code.  By the discussion preceding Theorem \ref{ratcan}, $\pr [b_R(v) \cap \mathcal{C} = \emptyset]$ is bounded above by
$$
\pr \left [ \left \{ \left ( \sum_{j=0}^{n-1} e_{j+1} \mathbf{x}^\intercal M_{\alpha,\mathcal{B}}^j \right )^{\!-1}\!\!\!\!w : \wt(w)=R \right \} \cap f_\mathcal{B}(v + \{\alpha^j\}_{j=1}^K) = \emptyset \right].
$$
The matrix in the left-hand term is uniformly distributed over all invertible matrices, by Lemma \ref{unif}.  Therefore, by Proposition \ref{invertible},
$$
\pr [b_R(v) \cap \mathcal{C} = \emptyset] \leq e^{-K} (c_q^{-1}+o(1)).
$$
Plugging this and $K = (R+1) \log n$ into (\ref{eq2}) yields
$$
\textrm{E}(|\Lambda_1(K)|+|\Lambda_2(K)|) \leq \frac{q^n \log n}{\binom{n}{R} (q-1)^R} (R+1 + o(1)),
$$
so a $(n,R,q)$-de Bruijn covering code of the desired length exists.
\end{proof}

\section{A Spectral Perspective} \label{section:spectral}

In this section, we describe a ``spectral'' test to see whether a given string is a de Bruijn covering code, and apply it to a probabilistic construction.  Define $e_N(x) = e^{2 \pi i x/N}$, as is standard notation.

\begin{prop} Let $S = (S(0),\ldots,S(M-1))$ be a $q$-ary string, for any $q > 1$.  Then $S$ is a de Bruijn covering code of radius $R$ and dimension $n$ if and only if the quantity
\begin{equation} \label{messeq}
\prod_{\omega = 0}^{q^n-1} \sum_{j=0}^{M-1} \sum_{v : \wt(v) \leq R} \sum_{m = 0}^{q^n - 1} e_{q^n}\left[ m ( \omega - \sum_{i=0}^{n-1} (S(i+j) + v_i \!\!\!\mod q) q^i ) \right]
\end{equation}
is positive, where $v$ varies over the set of $q$-ary sequences $(v_0,\ldots,v_{n-1})$ and the index of $S$ is written modulo $M$.  Otherwise, this expression is zero.
\end{prop}
\begin{proof} In what follows, all parameters vary over the ranges indicated in the statement above.  Note that
$$
\sum_m e_{q^n}(m(\omega - \omega^\prime))
$$
is positive if $\omega = \omega^\prime \!\! \mod q^n$, and zero otherwise.  If we represent a $q$-ary word as an integer base $q$, then the $j^{\textrm{th}}$ word appearing in $S$ is $\sum_i S(i+j) q^i$, and, if $\wt(v) \leq R$, this quantity plus $\sum_k v_k q^k$ (digits added independently modulo $q$) is the $j^\textrm{th}$ word with each symbol altered in at most $R$ coordinates.  Therefore, the quantity
$$
\sum_v \sum_m e_{q^n}(m ( \omega - \sum_i (S(i+j) + v_i \!\!\!\mod q) q^i ))
$$
is positive if and only if the word $S(i+j)$ is at most a distance $R$ from the word which is $\omega$ written base $q$.  Taking the sum over $j$ and then the product over $\omega$, we get that (\ref{messeq}) is positive if and only if $S$ is an $(n,R,q)$-de Bruijn covering code, and is zero otherwise.
\end{proof}

Consider the expected value of the above expression when we take a randomly, uniformly chosen binary string $S \in \{0,1\}^M$.  Clearly, an $(n,R,2)$-de Bruijn covering code of length $M$ exists if and only if this expected value is positive, since (\ref{messeq}) is always nonnegative.

\begin{theorem} An $(n,R,2)$-de Bruijn covering code of length $M$ exists if and only if
$$
\sum_{\textbf{j},\textbf{v},\textbf{m}} e_{2^n}\!\!\left[ \sum_{\omega=0}^{2^n-1} m_\omega (\omega - (2^n-1)/2) \right] \prod_{l=0}^{M-1} \cos\!\!\left( \pi \! \sum_{i,\omega} m_\omega (1 - 2v_{\omega,i}) 2^{i-n} \right) > 0,
$$
where $i$ and $\omega$ range over all pairs so that $0 \leq i \leq n-1$, $0 \leq \omega \leq 2^n-1$, and $i + j_\omega = l \!\!\!\mod M$, and the ranges of the other parameters are given by
$$ \begin{array}{c} \textbf{j} \in \{0,\ldots,M-1\}^{2^n} \\ \textbf{m} \in \{0,\ldots,2^n-1\}^{2^n} \\ \textbf{v} \in \{v \in \{0,1\}^n : \wt(v) \leq R\}^{2^n}.\end{array} $$
\end{theorem}
\begin{proof} First, rewrite (\ref{messeq}) by moving the product inside and collecting terms involving the same digits of $S$:
\begin{equation} \label{evenmessier}
\sum_{\textbf{j}} \sum_{\textbf{v}} \sum_{\textbf{m}} e_{2^n}\!\!\left[ \sum_{\omega=0}^{2^n-1} m_\omega \omega \right] \prod_{l=0}^{M-1} e_{2^n}\!\!\left[ - \sum_{i,\omega} m_\omega (S(l) + v_{\omega,i} \!\!\!\mod 2) 2^i \right].
\end{equation}
If $X$ is a random variable with two equally probable values $A$ and $B$, then $\textbf{E}[e_M(X)] = e_M((A+B)/2) \cos(\pi (A-B)/M)$.  Taking the expected value of (\ref{evenmessier}) therefore gives
$$
\sum_{\textbf{j},\textbf{v},\textbf{m}} e_{2^n}\!\!\left[ \sum_{\omega=0}^{2^n-1} m_\omega \omega \right] \prod_{l=0}^{M-1} e_{2^n}\!\!\left[ - \sum_{i,\omega} m_\omega 2^{i-1} \right] \!\cos\!\!\left( \pi \! \sum_{i,\omega} m_\omega (1 - 2v_{\omega,i}) 2^{i-n} \right)
$$
since the digits of $S$ are independent.  We may simplify this expression to
$$
\sum_{\textbf{j},\textbf{v},\textbf{m}} e_{2^n}\!\!\left[ \sum_{\omega=0}^{2^n-1} m_\omega (\omega - (2^n-1)/2) \right] \prod_{l=0}^{M-1} \cos\!\!\left( \pi \! \sum_{i,\omega} m_\omega (1 - 2v_{\omega,i}) 2^{i-n} \right).
$$
\end{proof}

Unfortunately, this result does not yield a practical means of calculating $M(n,R,2)$, due to the large number of terms.  Furthermore, it is unlikely that much cancellation can be identified in this sum, given the NP-hardness of determining a code's covering radius \cite{CHLL97}.  It may be possible, however, to exploit approximation algorithms for vertex-coverings to find a much simpler sum which yields a reasonable bound.

We also offer the following, in the spirit of the above results.

\begin{prop} Let $S = (S(0),\ldots,S(M-1))$ be a $q$-ary string, for any $q > 1$, and denote by $X$ the union of the radius $R$ balls about each codeword appearing as an $n$-string in $S$.  Then the number of points of $[q]^n$ not covered by $X$ is at most
$$
\sum_{\omega=0}^{q^n-1} \sum_{k=0}^{\infty} \frac{1}{k!} \left (-\sum_{j=0}^M \sum_{\wt(v) \leq R} \sum_{m=0}^{q^n-1} e_{q^n}\left[m (\omega - \sum_{i=0}^{n-1} (S(i+j) + v_{t,i} \!\!\!\mod q) q^i)\right]\right )^k
$$
where $v$ varies over the set of $q$-ary sequences $(v_0,\ldots,v_{n-1})$ and the index of $S$ is written modulo $M$.
\end{prop}
\begin{proof} As above, the quantity
$$
T(\omega) = q^{-n} \sum_{j=0}^M \sum_{\wt(v) \leq R} \sum_{m=0}^{q^n-1} e_{q^n}\left[m (\omega - \sum_{i=0}^{n-1} (S(i+j) + v_{t,i} \!\!\!\mod q) q^i)\right]
$$
counts the number of times that $\omega$ is covered.  Therefore $\sum_\omega e^{-q^n T(\omega)}$ is at least the number of uncovered points.
\end{proof}

One might conjecture that a sufficiently long sequence $S$ whose Fourier coefficients $\hat{S}(k)$ are small, for $k \neq 0$, covers all but a small fraction of Hamming space.  To avoid trivial cases, we must restrict our attention to sequences with approximately the same number of each symbol.  However, this statement is false even in the binary case, as illustrated by the following simple example.

Define $S = (S(0),\ldots,S(M-1))$, $M$ even, by $(S(2j),S(2j+1)) = (0,1)$ with probability $1/2$ and $(1,0)$ with probability $1/2$, each pair chosen independently.  Clearly, $S$ has the same number of $1$'s as $0$'s.  The $k^\textrm{th}$ Fourier coefficient, $k \neq 0$, has square magnitude
$$
|\hat{S}(k)|^2 = \sum_{u,v=0}^{M-1} e_M(k(u-v)) S(u) S(v).
$$
The values of $S(u)$ and $S(v)$ are independent if $|u-v|>1$, so the expected value of the above expression is
\begin{align*}
\textbf{E}[|\hat{S}(k)|^2] = & \sum_{u,v=0}^{M-1} e_M(k(u-v)) \textbf{E}[S(u) S(v)] \\
= & \sum_{u=0}^{M-1} \frac{1}{2} + \sum_{|u-v|>1} \frac{e_M(k(u-v))}{4} + \sum_{|u-v|=1} e_M(k(u-v)) \textbf{E}[S(u) S(v)] \\
\leq & \frac{M}{2} + \sum_{u,v=0}^{M-1} \frac{e_M(k(u-v))}{4} - \sum_{|u-v| \leq 1} \frac{e_M(k(u-v))}{4} + 2M \\
\leq & \frac{M}{2} + \left| \sum_{u=0}^{M-1} \frac{e_M(ku)}{2} \right|^2 + \frac{3M}{4} + 2M = \frac{15M}{4}.
\end{align*}
Any $n$-word appearing in $S$ has weight either $\lfloor n/2 \rfloor$ or $\lceil n/2 \rceil$.  Therefore, there exists a sequence $S$ of length $M$ with Fourier coefficients $\hat{S}(k) \ll \sqrt{M}$ so that, for any fixed $R$, the number of codewords at most a distance $R$ from the resulting code is an $O(n^{-1/2})$ fraction of the total.

It would be interesting to know whether the characteristic function of quadratic residues mod $p$ are a (near?) de Bruijn covering code whenever $p = \Omega(2^n/n^R)$.  Other possibilities for random-like constructions include the image of $[0,(p+1)/2]$ under the map $s \mapsto s^k$ with $(k,p-1)=1$, and the image of $[0,(p-1)/2]$ under the map $s \mapsto \tau^s$, for some primitive root $\tau$.  Unfortunately, because of the above example, the Fourier coefficients of these sets (which are known to be small) tell us nothing about how well they cover Hamming space.

\section{Numerical Bounds} \label{section:bounds}

It is of interest to know $M(n,R,q)$ for small values of its parameters -- in particular, for $q=2$, i.e., the binary case.  First, we collect a few simple observations.

\begin{enumerate}
\item $M(n,R,q) \leq M(n+k,R-l,q+m)$ for any $k,l,m \geq 0$.  If a de Bruijn covering code $\mathcal{C}$ exists for parameters $(n+k,R-l,q+m)$, then certainly decreasing the dimension, increasing the radius, or decreasing the number of symbols will leave $\mathcal{C}$ covering everything.  (In the case of decreasing the number of symbols, we can replace all occurrences of the excluded symbols to ``0''.  It is easy to check that this operation can only decrease distances from $n$-strings to the code.)
\item $M(n,0,q) = q^n$, as noted in the introduction.
\item $M(n,R,q) = 1$ if $R \geq n$, by taking the string ``0''.
\item $M(n,R,2) = 2$ if $\lfloor n/2 \rfloor \leq R < n$, by taking the string ``01''.  The two resulting codewords are complements in the $n$-cube, and therefore every string is within $\lfloor n/2 \rfloor$ of one of them.  Furthermore, it is clear that at least $2$ codewords are necessary.
\item $M(n,R,q) \geq K_q(n,R)$, the smallest number of codewords in a $q$-ary covering code of dimension $n$ and radius $R$.
\item $M(n,R,q) \neq M$ if $\min \{|n \mod M|,|(-n) \mod M|\} \leq n-2R-1$, where $|x \mod y|$ means the least nonnegative representative of $x$ modulo $y$.  Indeed, if a $(n,R,q)$-de Bruijn covering code $S=(s_0,\ldots,s_M)$ exists, then every string of $n$ consecutive symbols has weight
$$
\left \lfloor \frac{n}{M} \right \rfloor \wt(S) + \wt(s_i,\ldots,s_{i+A-1})
$$
for some $i$, where the indices are taken modulo $M$ and $A = |n \mod M|$.  Similarly, each such string has weight
$$
\left ( \left \lfloor \frac{n}{M} \right \rfloor + 1 \right) \wt(S) - \wt(s_i,\ldots,s_{i+B-1})
$$
for some $i$, where $B = |(-n) \mod M|$.  Therefore, any two codewords appearing in $S$ can differ by at most $C=\min \{A,B\}$ in weight.  If $C \leq n-2R-1$, then either the string $0^n$ or the string $1^n$ is at least a distance $R+1$ from any codeword.
\item Every $(n,R,2)$-de Bruijn covering code has a run of $\lfloor n/(R+1) \rfloor$ consecutive $0$'s and a run of $\lfloor n/(R+1) \rfloor$ consecutive $1$'s.  Suppose a code did not contain $0^k$ with $k = \lfloor n/(R+1) \rfloor$.  Then every element of the code has weight at least $\lfloor n / k \rfloor \geq R+1$, so the word $0^n$ is not covered, a contradiction.  An identical argument applies to the case of a run of $1$'s.
\item \label{possiblegaps} If there exists an $(n,R,q)$-de Bruijn covering code of length $M$, then there exists one of length $M+n+k-1$ for all $k \geq 0$.  If $S$ is the shorter string, append a copy of the first $(n-1)$ symbols and $k$ arbitrary $q$-ary symbols to the end.
\item \label{nogaps} If there exists an $(n,R,q)$-de Bruijn covering code of length $M(n,R,q)$ that somewhere contains the string $a^{n-1}$, then there exists an $(n,R,q)$-de Bruijn covering code of all lengths longer than $M(n,R,q)$.  We may simply insert more copies of $a$ into the string to generate longer ones.
\item There are at least $M(n,R,q)$ $(n,R,q)$-de Bruijn covering codes of length $M(n,R,q)$.  Since $M(n,R,q)$ is minimal, no such string has period less than $M(n,R,q)$, since otherwise we could truncate after a single period and achieve a smaller de Briujn covering code with the same parameters.  Therefore, all cyclic translations of any de Bruijn covering code -- which are each themselves de Bruijn covering codes -- are distinct.
\end{enumerate}

Below, we include a table of the best known bounds on the sizes of binary de Bruijn covering codes with various parameters.  A single number in an entry indicates that the exact value of $M(n,R,2)$ is known; two numbers indicate an upper and lower bound.  Bounds were achieved using the observations above, the table in \cite{Lit03}, as well as software that searched the string space randomly (for upper bounds), and one which searched it exhaustively (for lower bounds).  A few hundred hours of computing time on a 1.8 GHz Intel-based PC were used to construct this table.\\

\begin{table}[!ht] \label{tab:bounds}
\begin{center}
\begin{tabular}{|c||c|c|c|c|c|c|c|} \hline
$R \backslash n$ & 2 & 3 & 4 & 5 & 6 & 7 \\
\hline \hline 1  & 2 & 2 & 6 & 8 & 12 & 22 \\
\hline 2         & 1 & 2 & 2 & 2 & 8  & 10 \\
\hline 3         & 1 & 1 & 2 & 2 & 2  & 2 \\
\hline 4         & 1 & 1 & 1 & 2 & 2  & 2 \\
\hline 5         & 1 & 1 & 1 & 1 & 2  & 2 \\
\hline 6         & 1 & 1 & 1 & 1 & 1  & 2 \\
\hline \multicolumn{7}{c}{} \\\hline $R \backslash n$ & 8 & 9 & 10 & 11 & 12 & 13 \\
\hline \hline 1  & 32    & 57-130 & 105-322 & 180-694 & 342-1454 & 598-2937 \\
\hline 2         & 14    & 20     & 38      & 38-117  & 62-244   & 97-529 \\
\hline 3         & 6     & 12     & 16      & 20      & 34-40    & 34-119 \\
\hline 4         & 2     & 2      & 4       & 8       & 16       & 24  \\
\hline 5         & 2     & 2      & 2       & 2       & 8        & 8 \\
\hline 6         & 2     & 2      & 2       & 2       & 2        & 2  \\
\hline 7         & 2     & 2      & 2       & 2       & 2        & 2 \\
\hline 8         & 1     & 2      & 2       & 2       & 2        & 2 \\
\hline 9         & 1     & 1      & 2       & 2       & 2        & 2 \\
\hline 10        & 1     & 1      & 1       & 2       & 2        & 2 \\
\hline 11        & 1     & 1      & 1       & 1       & 2        & 2 \\ \hline
\end{tabular}\caption{Best known bounds for $M(n,R,2)$}
\end{center}
\end{table}

\section{Remarks and Further Questions} \label{section:remarksandfurtherquestions}

Statement \ref{possiblegaps} in the previous section highlights a frustrating property of de Bruijn covering codes that stands in stark contrast to ordinary covering codes: it is possible for one to exist of length $M$ but for none to exist of length $M+1$.  For example, a $(10,4,2)$ code exists of lengths $4$ (``1100''), $6$ (``011100''), $8$ (``00111100''), and $12$ (``000011111100''), but none of lengths $5$, $7$, $9$, $10$, or $11$ exist.  However, by the above, a $(10,4,2)$ code of all lengths at least $13$ must exist.  Therefore, in addition to finding the smallest possible de Bruijn covering code, we would like to know when de Bruijn covering codes with lengths \textit{between} $M(n,R,q)$ and $M(n,R,q)+n-1$ exist.

Another difference between de Bruijn covering codes and ordinary ones is that there is no easy way to use known  efficient codes to build efficient codes for larger $n$, smaller $R$, or larger $q$.  It would be desirable to define a ``product'' analogous to direct sums for ordinary covering codes.  Unfortunately, interlacing, the obvious candidate for such a product, appears to be very inefficient.  We offer a different, though related construction which allows us to increase $q$ when the desired number of symbols is a perfect power of the number of symbols in the original code.

\begin{prop} \label{enlarge} If $a^s = b$ for any positive integers $a$, $b$, and $s$, then for all $n,R > 0$,
$$
M(n,R,b) \leq s^2 \left \lceil \frac{M(sn,R,a)+sn}{s} \right \rceil - s.
$$
\end{prop}
\begin{proof} Let $t = M(sn,R,a)$ and $m = s^2 \lceil (t+sn)/s \rceil - s$, and let $C = (c_0,\ldots,c_{t-1})$ be a minimum-length $(sn,R,a)$-de Bruijn covering code.  We construct an $(n,R,a^s)$-de Bruijn covering code $C^\prime = (c^\prime_0,\ldots,c^\prime_{m-1})$ of length $m$.  Choose some bijection $\sigma$ between $(\mathbb{Z}/a \mathbb{Z})^s$ and $\mathbb{Z}/a^s \mathbb{Z}$, and define
$$
c^\prime_j = \sigma(c_{|sj \!\!\! \mod (m/s)|},\ldots,c_{|s(j+1)-1 \!\!\! \mod (m/s)|})
$$
with indices on the left hand side taken modulo $m$ and indices on the right hand side taken modulo $t$.  Evidently, $C^\prime$ is well defined, since $s | m$.  Now, suppose $X=(x_0,\ldots,x_{n-1})$ is an $n$-string over $a^s$ symbols.  We claim that there is some codeword in the set of consecutive $n$-strings of $C^\prime$ which is within $R$ symbols of $x$.

Indeed, let $x^\prime_j = \sigma^{-1}(x_j)$ for $0 \leq j < n$ and define $X^\prime = x^\prime_0 \cdots x^\prime_{sn-1}$, a string of length $sn$.  Then some string $X^{\prime\prime}$ which differs from $X^\prime$ in at most $R$ symbols occurs somewhere in $C$, say, beginning at coordinate $k$.  $X^{\prime\prime}$ must occur at least $s$ times in $C^\prime$, at coordinates $k+jm/s$ for $0 \leq j < s$.  (If $X^{\prime\prime}$ ``wraps around'' in $C$, the extra $\geq sn-1$ symbols at the end of each block of length $m/s$ guarantee $X^{\prime\prime}$ appears in $C^\prime$.)  Furthermore, since $(m/s,s)=1$, the numbers $k + jm/s$, $0 \leq j < s$, represent all residue classes modulo $s$, so there is some $r$ so that $k + rm/s \equiv 0 \mod s$.  Then the string
$$
\sigma^{-1}(c^\prime_{k+rm/s},\ldots,c^\prime_{k+rm/s+s-1})\ldots\sigma^{-1}(c^\prime_{k+rm/s+(n-1)s},\ldots,c^\prime_{k+rm/s+ns-1}) $$
appears in $C$ and at most $R$ of its coordinates differ from those of $X$.
\end{proof}

The most obvious question arising from the subject of the present work is the issue of whether the bound stated in Theorem $\ref{mainthm}$ is best possible, i.e., whether the $\log$ factor can be dropped or the result can be extended to $q$'s which are not prime powers.  We also would like to explain why so many of the entries in Table 1 are even.


\begin{thebibliography}{100}
\bibitem {AS00} N.\@ Alon, J.\@ Spencer, The probabilistic method. Wiley-Interscience Series in Discrete Mathematics and Optimization. Wiley-Interscience [John Wiley \& Sons], New York, 2000.
\bibitem {Beard86} T.\@ B. \@ Beard, Jr.\@, Matrix fields, regular and irregular: a complete fundamental characterization. Linear Algebra Appl. {\bf 81} (1986), 137--152.
\bibitem {CEK02} J.\@ N.\@ Cooper, R.\@ B.\@ Ellis, A.\@ B.\@ Kahng, Asymmetric binary covering codes.  J. Combin. Theory Ser. A {\bf 100} (2002), no. 2, 232--249.
\bibitem {CHLL97} G.\@ Cohen, I.\@ Honkala, S.\@ Litsyn and A.\@ Lobstein, Covering codes. North-Holland Mathematical Library 54, Elsevier, 1997.
\bibitem {Fred82} H.\@ Fredricksen, A survey of full length nonlinear shift register cycle algorithms. SIAM Rev. 24 (1982), no. 2, 195--221.
\bibitem {HJ00} R.\@ A.\@ Horn, C.\@ R.\@ Johnson, Matrix analysis. Cambridge University Press, Cambridge, 1990.
\bibitem {Hoch95} D.\@ Hochbaum, ed., Approximation Algorithms for NP-Hard Problems, PWS Publishing Company, Boston, MA, 1995.
\bibitem {KSV03} M.\@ Krivelevich, B.\@ Sudakov, V.\@ Vu, Covering codes with improved density. Preprint, 2003.
\bibitem {Land00} M.\@ Landsberg, Feedback functions for generating cycles over a finite alphabet. Discrete Math. {\bf 219} (2000), no. 1-3, 187--194.
\bibitem {Lit03} S.\ Litsyn, Table of the best currently known lower and upper bounds on the smallest size of a covering code, Manuscript, {\bf http://www.eng.tau.ac.il/$\sim$litsyn/tablecr/index.html}.
\bibitem {Park50} W.\@ V.\@ Parker, The matrix equation $AX=XB$. Duke Math. J. {\bf 17} (1950), 43--51.
\end{thebibliography}
\end{document}